\documentstyle[12pt]{article}

\newtheorem{prop}{Proposition}[section]
\newtheorem{dfn}[prop]{Definition}
\newtheorem{theo}[prop]{Theorem}

\newtheorem{conj}[prop]{Conjecture}

\newtheorem{coro}[prop]{Corollary}

\def\R{{\bf R }}

\def\Z{{\bf Z }}
\def\Q{{\bf Q }}
\def\P{{\bf P }}

\def\g{ \gamma}

\def\G{ \Gamma }
\def\p{ \varphi }

\def\ra{\rightarrow}

\title{Rational points on some Fano cubic bundles}

\author{
Victor V. Batyrev\thanks{Supported by Deutsche Forschungsgemeinschaft} \\
\small  Universit\"at-GHS-Essen, Fachbereich  6,  Mathematik \\
\small  Universit\"atsstr. 3,  45141  Essen, Germany  \\
\small  e-mail: victor.batyrev@uni-essen.de \\
and \\
Yuri Tschinkel\thanks{Leibniz Fellow of the EC at ENS, Paris}\\
\small Dept. of Mathematics, U.I.C.\\
\small Chicago, (IL) 60608,  U.S.A.  \\
\small e-mail: yuri@math.uic.edu
}

\begin{document}

\date{}

\maketitle

\thispagestyle{empty}

\begin{abstract}
We consider
smooth Fano hypersurfaces  $X_{n+2} \subset \P^n \times {\P}^3 $ $(n \geq
1)$ given by a polynomial
$$\sum_{i =0}^3 l_i({\bf x})y_i^3 \in
\Q [ x_0, \ldots, x_n, y_0,\ldots, y_3 ]$$
where $l_0({\bf x}), \ldots, l_3({\bf x})$ are
homogeneous linear forms in $x_0, \ldots, x_n$.
We obtain lower bounds
for the number of $F$-rational points  of bounded
anticanonical height in arbitrary nonempty Zariski open subsets
$U \subset X_{n+2}$ for number fields $F$ containing
$\Q(\sqrt{-3})$.  These bounds contradict
previous expectations about the distribution of $F$-rational points
of bounded height on Fano varieties.
\end{abstract}


\section{Cubic bundles}

Let $X_{n+2}$ be a hypersurface in $\P^n \times \P^3$ $(n \geq
1)$ defined by the equation
$$ P({\bf x},{\bf y}) = \sum_{i =0}^3 l_i({\bf x})y_i^3 = 0$$
where
$$P({\bf x},{\bf y}) \in \Q [ x_0, \ldots, x_n,
y_0,\ldots, y_3 ]$$
and $l_0({\bf x}), \ldots, l_3({\bf x})$ are
homogeneous linear forms in $x_0, \ldots, x_n$.
Put $k = \max ( n+1, 4)$.
We shall always assume that any $k$ forms among
 $ l_0({\bf x}), \ldots, l_3({\bf x})$ are linearly independent.
It is elementary to
check the  following statements:

\begin{prop}
The hypersurface $X_{n+2}$ is a smooth Fano variety
containing  a Zariski open subset
$U_{n+2}$ which is isomorphic to  ${\bf A}^{n+2}$.
\end{prop}

\begin{prop}
Let  $U_P \subset \P^n$ be the Zariski open subset
defined by the condition
\[ \prod_{i =0}^3 l_0({\bf x}) \neq 0. \]
Then the  fibers of the natural projection  $\pi \,:\, X_{n+2} \ra \P^n$
 over closed points of $U_P$ are
smooth diagonal cubic surfaces in $\P^3$.
\end{prop}

\noindent
By Lefschetz theorem, we have:

\begin{prop}
The Picard group of $X_{n+2}$ over an arbitrary field containing $\Q$ is
isomorphic to ${\Z}\oplus \Z$.
\label{pic}
\end{prop}

\section{Heights on cubic surfaces}

Let $F$ be a number field, ${\rm Val}(F)$ the set of all
valuations of $F$, $W$ a projective algebraic variety over $F$,  $W(F)$
the set of $F$-rational points of $W$,
$D$ a very ample divisor on $W$, and ${\g} =
\{ s_0, \ldots, s_m\}$ a basis over $F$
of the space of global sections $\G(W, {\cal O}(D))$.
The {\em height function associated with $D$ and
${\g}$}
$$H(W, D,{\g},x)\; : \; W(F) \ra \R_{>0} $$
is given by the formula
\[ H(W, D,{\g},x) = \prod_{v \in {\rm Val}(F)}
\max_{i=0, \ldots, m} | s_i(x) |_v,  \]
where $|\cdot |_v \,:\, F_v \ra {\R}_{>0}$ is the multiplier
of a   Haar measure on the additive group of the
$v$-adic completion of $F$.

\begin{dfn}
{\rm Let $Z \subset W$ be a locally closed algebraic
subset of $W$, $B$ a positive
real number. Define
\[ N(Z, D,{\g},B) := {\rm Card}\{ x \in W(F) \cap Z \mid
 H(W, D,{\g},x) \leq B \}.\]}
\end{dfn}

\noindent
The following classical statement is due to A. Weil:

\begin{theo}
Let ${\g'} = \{ s_0', \ldots, s_m' \}$ be another basis in
$\G(W, {\cal O}(D))$.
Then there exist two positive constants $c_1,c_2$ such that
\[ c_1 \leq  \frac{H(W, D,{\g},x)}{H(W, D,{\g'},x)} \leq  c_2  \]
for all $x \in W(F)$.
\label{weil}
\end{theo}

For a smooth projective variety $W$ we denote by  $-K_W$
the anticanonical divisor on $W$.

\begin{theo}
Let $Y \subset {\bf P}^3$ be a smooth cubic surface defined over $F$,
${\g} = \{ s_0, s_1, s_2, s_3 \}$ the basis of
global sections of ${\cal O}(-K_Y)$ corresponding to the standard
homogeneous coordinates on ${\P}^3$.
Assume that $Y$ can be obtained by blowing up of $6$ $F$-rational
points in $\P^2$.  Then for any nonempty
Zariski open subset $U \subset Y$ one has
\[ N(U, -K_Y, {\g}, B) \geq c B (\log B)^3  \]
for all $B > 0$ and for some positive constant $c$.
\label{cub}
\end{theo}

\noindent
{\em Proof.} Assume that $Y$ is obtained by blowing up of
$p_1, \ldots, p_6 \in \P^2(F)$.  By \ref{weil},
we can assume  without loss
of generality that $p_1 = (1:0:0)$, $p_2 = (0:1:0)$ and
$p_3 = (0:0:1)$. Denote by
$Y_0$ the Del Pezzo surface obtained by blowing up $p_1, p_2, p_3$.
Let $f \,: \, Y \ra Y_0$ be the contraction
of exceptional curves $C_4, C_5, C_6 \subset Y$ lying over $p_4,p_5,p_6$.

Let $V$ be the $10$-dimensional  space over $F$ of all homogeneous
polynomials of degree
$3$ in variables $z_0, z_1, z_2$. We  identify
$\G(Y, {\cal O}(-K_Y))$ with the subspace in $V$ consisting
of all polynomials vanishing in $p_1, \ldots, p_6$.
Analogously, we  identify
$\G(Y_0, {\cal O}(-K_{Y_0}))$ with the subspace in $V$ consisting
of all polynomials vanishing in $p_1,p_2,p_3$.  Let
$\g_0 = \{ s_0, \ldots, s_6 \} \subset V$ the extension of the basis
${\g}$ to a basis of the subspace  $\G(Y_0, {\cal O}(-K_{Y_0})) \subset V$.

The surface  $Y_0$ is an smooth equivariant compactification of the split
$2$-dimensional algebraic torus over $F$
$$({\bf G}_{m})^2 = \P^2 \setminus
\{ l_{12}, l_{13}, l_{23} \}$$
where $l_{ij}$ denotes the projective line in $\P^2$ through $p_i$ and $p_j$.
Since $Y_0$ is a smooth toric variety,
the main theorem in \cite{BaTschi1} shows that
the following asymptotic formula holds:
\begin{equation}
N(({\bf G}_{m})^2, -K_{Y_0}, \g_0', B) =
c_0 B (\log B)^3(1 + o(1)),\;\;
B \ra \infty,
\label{f0}
\end{equation}
where $c_0$ is some positive constant and
\[ \g_0' = \{ z_0z_1z_2, z_1^2z_2, z_1 z_2^2, z_2^2z_0, z_2z_0^2,
z_0^2z_1, z_0 z_1^2 \}. \]

Let $U$ be any nonempty Zariski open subset in $Y$. We denote by
$U_0$ a nonempty open subset in $U$ such that the restriction
of $f$ on $U_0$ is an isomorphism and $f(U_0)$ is contained in
$({\bf G}_{m})^2 \subset Y_0$. Since
\[  \prod_{v \in {\rm Val}(F)}
\max_{i=0, \ldots, 3} | s_i(x) |_v  \leq
\prod_{v \in {\rm Val}(F)}
\max_{i=0, \ldots, 6} | s_i(x) |_v  \]
holds for every $F$-rational point $x \in U_0$, we
obtain
\begin{equation}
N(U_0, -K_Y, \g, B) \geq N(U_0, -K_{Y_0}, \g_0, B)
\label{f1}
\end{equation}
for any $B >0$. By \ref{weil}, there exists a positive constant
$c_3$ such that
\begin{equation}
 N(U_0, -K_{Y_0}, \g_0, B) \geq c_3 N(U_0, -K_{Y_0}, \g_0',
B).
\label{f2}
\end{equation}
On the other hand,
\begin{equation}
 N(U_0, -K_{Y_0}, \g_0',
B)  = N(({\bf G}_{m})^2, -K_{Y_0}, \g_0', B) - N(Z, -K_{Y_0}, \g_0',
B),
\label{f3}
\end{equation}
where $Z = ({\bf G}_{m})^2 \setminus U_0$. Let $Z_1, \ldots,
Z_l$ be the irreducible components of $Z$, $\overline{Z}_i$
the closure of $Z_i$ in $Y_0$ $(i =1, \ldots, l)$. It is known
that
\begin{equation}
 N(Z_i, -K_{Y_0}, \g_0',B) \leq  c_4 B^{2/({\rm deg}\,
\overline{Z}_i)}
\label{f4}
\end{equation}
holds for some positive constant $c_4$,
where ${\rm deg}\,\overline{Z}_i$ denotes the degree of
$\overline{Z}_i$ with respect to the anticanonical divisor
$-K_{Y_0}$. Since every irreducible curve  $C \subset Y_0$ with
${\rm deg}\, C = 1$ is a component of $Y_0 \setminus ({\bf G}_m)^2$, we
have ${\rm deg}\,\overline{Z}_i \geq 2$; i.e.,
\begin{equation}
 N(Z_i, -K_{Y_0}, \g_0',B) \leq  c_4 B
\label{f5}
\end{equation}
holds  for all $i = 1, \ldots, l$.

It follows from the asymptotic formula (\ref{f0}) combined with
 (\ref{f1}),  (\ref{f2}) and (\ref{f3})
that there exists a positive constant $c$ such that
\[ N(U_0, -K_Y, {\g}, B) \geq c B (\log B)^3  \]
holds for all $B > 0$. This yields the statement, since $U_0$ is
contained in $U$.
\hfill $\Box$

\begin{coro}
Let
$Y$ be a smooth diagonal
cubic surface in $\P^3$ defined by the equation
\[ a_0 y_0^3 + a_1 y_1^3 + a_2 y_2^3 + a_3 y_3^3 = 0 \]
with coefficients $a_0, \ldots, a_3$ in  a number field
$F$ which contains ${\Q}(\sqrt{-3})$. Assume that
there exist numbers $b_0, \ldots, b_3 \in F^*$ such that
$a_i = b_i^3$ $(i =0, \ldots, 3)$. Then for any nonempty
Zariski open subset $U \subset Y$ one has
\[ N(U, -K_Y,  \g, B) \geq c B (\log B)^3  \]
for all $B > 0$ and some positive constant $c$.
\label{cubic1}
\end{coro}

\noindent
{\em Proof.} It follows from our assumptions on the coefficients
$a_0, \ldots, a_3$ and on the field $F$ that
all $27$ lines on $Y$ are  defined over $F$. Hence,
$Y$ can be obtained from ${\P}^2$ by blowing up of $6$
$F$-rational points.  Now the statement follows from
\ref{cub}. \hfill $\Box$

\section{Rational points on $X_{n+2}$}

We have the natural isomorphism
\[ \G (X_{n+2}, {\cal O}(-K_{X_{n+2}}) ) \cong
\G (\P^3, {\cal O}(1)) \otimes \G (\P^n, {\cal O}(n)).\]
Let $\{ t_0, \ldots, t_m \}$ be a basis
in $\G (\P^n, {\cal O}(n))$ and  $\{ s_0, s_1, s_2, s_3 \}$ the
standard basis in $\G (\P^3, {\cal O}(1))$. Denote by
$\g$ the basis of $\G (X_{n+2}, {\cal O}(-K_{X_{n+2}}) )$
consisting of $s_i \otimes t_j$ $(i=0,\ldots, 3; \; j =
0, \ldots, m)$.

\begin{theo} Let $n \geq 3$. Then
for any nonempty Zariski open subset $U \subset X_{n+2}$ and for
any field $F$ containing  ${\Q}(\sqrt{-3})$, one has
\[ N(U, -K_{X_{n+2}},  \g, B) \geq c B (\log B)^3  \]
for all $B > 0$ and some positive constant $c$.
\label{t1}
\end{theo}

\noindent
{\em Proof.} Consider the projection $\pi\,:\,U\ra \P^n$.
Let $U'$ be a Zariski open subset in $\pi (U) \cap U_P$  such that
the image of $U'$ under the dominant (rational)  mapping
\[ \psi \;:\; \P^n - - \ra \P^3  \]
\[ \psi(x_0 : \ldots : x_n) =
(l_0({\bf x}): \ldots : l_3({\bf x})) \]
is Zariski open in $\P^3$.
Denote by
$\p$ the finite morphism
\[ \p \; : \; \P^3 \ra \P^3, \]
\[ \p(z_0 : \ldots : z_3) = (z_0^3 : \ldots : z_3^3). \]
Since $\P^3(F)$ is Zariski dense in $\P^3$,
there is a point $p \in \P^3(F) \cap \p^{-1}(\psi(U'))$.
Since $U'(F) \cap \psi^{-1}(\p(p))$
is Zariski dense  in $\psi^{-1}(\p(p))$,
there exists  $q \in U'(F) \cap
\psi^{-1}(\p(p))$. Therefore, the  fiber of $\pi $ over $q$
is a diagonal cubic surface $Y_q$ such that and $U \cap Y_q \subset Y_q$ is
a nonempty Zariski open subset. It remains to apply \ref{cubic1}.
\hfill $\Box$

\begin{theo}
Let  $n =2$. Then  there exists a
number field $F_0$ depending only on $X_{n+2}$ such that
for any nonempty Zariski open subset $U \subset X_{n+2}$ for any field
$F$ containing $F_0$ one has
\[ N(U, -K_{X_{n+2}},  \g, B) \geq c B (\log B)^3  \]
for all $B > 0$ and some positive constant $c$.
\label{t2}
\end{theo}

\noindent
{\em Proof.}
Let $U'$ be a Zariski open subset in $\pi (U) \cap U_P$. We have the
linear embedding
\[ \psi \;:\; \P^2 \hookrightarrow   \P^3 \]
defined by $l_0({\bf x}), \ldots, l_3({\bf x})$.
Then $\p^{-1}(\psi(\P^2))$ is a smooth diagonal cubic surface
$S \subset  \P^3$ defined over $\Q$. Let $F_0$ be a finite
extension of $\Q(\sqrt{-3})$ such that $S(F_0)$ is Zariski dense in $S$.
Then there exists a point $p \in S(F_0)$ such that $q = \p(p)$ is contained
in $U'$. Therefore,
the fiber of $\pi$ over $q$
is a diagonal cubic surface $Y_q$, and $U \cap Y_q \subset Y_q$ is
a nonempty Zariski open subset. It remains to apply
\ref{cubic1}.
\hfill $\Box$

\begin{theo}
Let  $n = 1$. Then for
any nonempty Zariski open subset $U \subset X_{n+2}$, there exists a
number field $F_0$ $($which depends on $U$$)$ such that for any field
$F$ containing $F_0$, one has
\[ N(U, -K_{X_{n+2}},  \g, B) \geq c B (\log B)^3  \]
for all $B > 0$ and some positive constant $c$.
\label{t3}
\end{theo}

\noindent
{\em Proof.}
Let $U'$ be a Zariski open subset in $\pi (U) \cap U_P$. We have the
linear embedding
\[ \psi \;:\; \P^1 \hookrightarrow  \P^3 \]
defined by $l_0({\bf x}), \ldots, l_3({\bf x})$.
Then $\p^{-1}(\psi(\P^1))$ is an algebraic curve
$C \subset  \P^3$ which is a complete intesection of two
diagonal cubic surfaces defined over $\Q$.
Let $F_0$ be a finite extension of $\Q(\sqrt{-3})$ such that
there exists an $F_0$-rational point  $p \in C(F_0) \cap
\p^{-1}(U')$.
Then the fiber of $\pi$ over $q = \p(p)$
is a diagonal cubic surface $Y_q$ and $U \cap Y_q \subset Y_q$ is
a nonempty Zariski open subset. It remains to apply
\ref{cubic1}.
\hfill $\Box$

\section{Conclusions}

The following statement, inspired by the Linear Growth conjecture
of Manin (\cite{manin1}) and by extrapolation
of known results (circle method, flag varieties,
toric varieties),  has been expected to be true
\cite{bat.man,franke-manin-tschinkel}:

\begin{conj}
Let $X$ be a smooth Fano variety over a number field $E$.
Then there exist
a Zariski open subset $U \subset X$ and a finite
extension $F_0$ of $E$ such that for all number fields $F$
containing $F_0$ the following
asymptotic formula holds
\[ N(U, -K_X, \g, B) = c B (\log B)^{t-1} (1 + o(1)), \;\;
B \ra \infty ,\]
where $t$ is the rank of the Picard group of
$X$ over $F$.
\label{conjecture}
\end{conj}

Some lower and upper bounds for $N(U, -K_X, \g, B)$ for Del Pezzo
surfaces and Fano threefolds have been obtained in
\cite{manin1,MaTschi}.

The  conjecture \ref{conjecture} was refined by E. Peyre who proposed
an adelic interpretation for the constant $c$ introducing  Tamagawa numbers
of Fano varieties \cite{peyre}. This refined version of
the conjecture has been proved for toric varieties in
\cite{BaTschi,BaTschi1}.

The statements in Theorems \ref{t1}, \ref{t2}, \ref{t3} and the
property \ref{pic} show
that Conjecture \ref{conjecture} is  not true for Fano
cubic bundles $X_{n+2}$ $(n \geq 1)$.

\end{document}